\documentclass{birkjour}
\usepackage[cp1251]{inputenc}  
\usepackage[T2A]{fontenc}      
\usepackage{amssymb}
\usepackage{amssymb,amsmath}
\newtheorem{thm}{Theorem}[section]
\newtheorem{lem}[thm]{Lemma}
\newtheorem{cor}[thm]{Corollary}

\begin{document}

\title{\textbf{On the Skitovich-Darmois theorem for ${\boldsymbol a}$-adic solenoids}}
\author{I.P. Mazur}

\address{
47 Lenin Ave.,\\
Kharkov 61103,\\
Ukraine}

\email{mazurivan85@gmail.com}

\keywords{the Skitovich-Darmois theorem, a functional equation, an
${\boldsymbol a}$-adic solenoid}

 \maketitle

\begin{abstract}

Let $X$ be a compact connected Abelian group. It is well-known that then
there exist topological automorphisms
  $\alpha_j, \beta_j $ of $X$ and independent random variables $\xi_1$
and $\xi_2$ with values
  in $X$ and distributions $\mu_1, \mu_2$
   such that
the linear forms $L_1 = \alpha_1\xi_1 + \alpha_2\xi_2$ and $L_2 =
\beta_1\xi_1 + \beta_2\xi_2$ are independent, whereas $\mu_1$ and
$\mu_2$ are not represented as convolutions of Gaussian and
idempotent distributions. This means that  the Skitovich--Darmois
theorem fails for such groups. We prove that if we consider three
linear forms of three independent random variables taking values
in $X$, where $X$ is an ${\boldsymbol a}$-adic solenoid, then the
independence of the linear forms implies that at least one of the
distributions is idempotent. We describe all such solenoids.

\end{abstract}

\section{Introduction}

It is well-known that proofs of many characterization theorems of
mathematical statistics are reduced to solving of some functional
equations. Consider the classical Skitovich--Darmois theorem that
characterizes Gaussian distributions on the real line (\cite[ch.
3]{KLR}): Let $\xi_i ,i=1,2,\ldots,n,$ $n\geq2,$ be independent
random variables, and $\alpha_j,\beta_j$ be nonzero constants.
Suppose that the linear forms
$L_1=\alpha_1\xi_1+\cdots+\alpha_n\xi_n$ and
$L_2=\beta_1\xi_1+\cdots+\beta_n\xi_n$ are independent. Then all
random variables $\xi_j $ are Gaussian.

Let $\hat{\mu}_j(y)$ be the characteristic functions of the
distributions of   $\xi_j,j=1,2,\ldots,n.$ Taking in the account
that $\mathbf{E}[e^{i\xi_jy}]=\hat{\mu}_j(y)$, it is easy to
verify that the Skitovich--Darmois theorem is equivalent to the
following statement: The solutions of the Skitovich--Darmois
equation

$$\prod_{j=1}^{n} \hat{\mu}_j(\alpha_j u+\beta_j v)=\prod_{j=1}^{n}
\hat{\mu}_j(\alpha_j u)\hat{\mu}_j(\beta_j v),\quad u,v\in
\mathbb{R},$$ in the class of the normalized continuous positive
definite functions are the characteristic functions of the
Gaussian distributions, i.e. $\hat{\mu}_j(y)=\exp\{ia_jy-\sigma_j
y^2\}, \ a_j\in\mathbb{R},\sigma_j \geq 0,y\in\mathbb{R}, \
j=1,2,\ldots,n$.

This theorem was generalized to various classes of locally compact
Abelian groups (see for example \cite{F 1992}--\cite{FM 2008},
\cite{Mazur Fin}). In these researches random variables take
values in a locally compact Abelian group $X$, and coefficients of
the linear forms are topological automorphisms of $X$. As in the
classical case the characterization problem is reduced to the
solving of the Skitovich--Darmois equation in the class of the
normalized continuous positive definite functions on the character
group of the group $X$.

In \cite{FG 2000} Feldman and Graczyk have shown, that  even a
weak analogue of the Skitovich--Darmois theorem fails for compact
connected Abelian groups. Namely, they proved the following
statement: Let $X$ be an arbitrary compact connected Abelian
group. Then there exist topological automorphisms
$\alpha_j,\beta_j,j=1,2,$ of $X$ and independent random variables
$\xi_1, \xi_2$ with values in $X$ and having  distributions, that
are not convolutions of the Gaussian and idempotent distributions,
wereas the linear forms $L_1=\alpha_1\xi_1+\alpha_2\xi_2$ and
$L_2=\beta_1\xi_1+\beta_2\xi_2$ are independent.

The aim of this article is to show that a weak analogue of the
Skitovich--Darmois theorem holds for some compact connected
Abelian groups if we consider three linear forms of three random
variables. Namely, we will construct an ${\boldsymbol a}$-adic
solenoid $\Sigma_{\boldsymbol a}$ (we give the  full description
of such solenoids in the Theorem \ref{Theorem Sol p}) for which
the independence of three linear forms of three independent random
variables with values in $\Sigma_{\boldsymbol a}$ implies that at
least one random variable has idempotent distribution.

\section{Definitions and notation}

Let $X$ be a second countable locally compact Abelian group.
Denote by $Aut(X)$ the group of the topological automorphisms of $X$.
Let $k$ be an integer. Denote by $f_k$ the mapping $f_k:X\rightarrow X$
definite by the equality $f_kx=kx$. Put $X^{(k)}=f_k(X)$.

Let $Y=X^*$ be the character group of $X$. The value of a
character $y\in Y$ at $x\in X$ denote by $(x,y)$. Let $B$ be a
nonempty subset of $X$. Put
$$A(Y,B)=\{y\in Y:(x,y)=1 , x\in\ B\}.$$
The set $A(Y,B)$ is called the annihilator of $B$ in $Y$. The
annihilator $A(Y,B)$ ia a closed subgroup in $Y$. For each
$\alpha\in Aut(X)$ definite the mapping $\tilde{\alpha}: Y
\rightarrow Y$ by the equality $(\alpha x,y)=(x,\tilde{\alpha}y)$
for all $x\in X, y \in Y$. The mapping $\tilde{\alpha}$ is a
topological automorphism of $Y$. It is called an adjoint of
$\alpha$. The identity automorphism  of a group $X$
denote by $I$.

In the paper we will use standard facts of abstract harmonic
analysis  (see \cite{H R}). Let $\mu$ be a distribution on $X$.
The characteristic function of $\mu$ is definite by the formula
$$\hat{\mu}(y)=\int_X(x,y)d\mu(y),y\in Y.$$
Put $F_{\mu}=\{y\in Y: \hat{\mu}(y)=1\}$. Then $F_{\mu}$ is a
subgroup of $Y$, and the function $\hat{\mu}(y)$ is $F_{\mu}$ -
invariant, i.e. $\hat{\mu}(y+h)=\hat{\mu}(y),y\in Y,h\in F_{\mu}$.

Denote by $E_x$ the degenerate distribution  concentrated at $x$.
Let $K$ be a compact subgroup of $X$. Denote by $m_K$ the Haar
distribution on $K$. Denote by $I(X)$ the set of shifts of such
distributions,  i.e. the distributions of the form $m_K\ast E_x$,
where $K$ is a compact subgroup of $X$, $x\in X$. The
distributions of the class $I(X)$ are called idempotent. Note that
the characteristic function of $m_K$ is of the form:
$$
\hat{m}_K(y)=
\begin{cases}
1,& y\in A(Y,K), \\
0, & y\not\in A(Y,K).
\end{cases}
$$

A distribution $\mu$ on the group $X$ is called Gaussian
(\cite{Parth}) if its characteristic function can be represented
in the form
$$\hat{\mu}(y)=(x,y)\exp\{-\varphi(y)\}, \quad y\in Y,$$
where $\varphi(y)$ is a continuous nonnegative function satisfying
the equation
$$\varphi(u+v)+\varphi(u-v)=2(\varphi(u)+\varphi(v)),\quad u,v\in
Y.$$  Denote by $\Gamma(X)$ the set of Gaussian distributions on
$X$.

Denote by $\mathbb{Z}$ the
infinite cyclic group, by $\mathbb{R}$ the additive group of real
numbers, by $\mathbb{T}$ the circle group, by $\mathbb{Q}$ the additive
group of rational numbers with the discrete topology, by $\Delta_{{\boldsymbol a}}$ the group
of ${\boldsymbol a}$-adic integers, by
$\mathbb{Z}(m)$ the group of residue modulo $m$.

Let ${\boldsymbol a}=(a_0,a_1,\ldots,a_n,\ldots)$ be a fixed but
arbitrary infinite sequence of natural numbers, where all $a_i>1$.
Consider the group $\mathbb{R}\times \Delta_{{\boldsymbol a}}$.
Let $B$ be a subgroup of $\mathbb{R}\times \Delta_{{\boldsymbol
a}}$ of the form $B=\{(n,n\mathbf{u})\}_{n=-\infty}^{\infty}$,
where $\mathbf{u}=(1,0,\ldots,0,\ldots)$. The factor-group
$\mathbf{\Sigma}_{\boldsymbol a}=(\mathbb{R}\times
\Delta_{\boldsymbol a})/B$ is called an ${\boldsymbol a}$-adic
solenoid. The group $\mathbf{\Sigma}_{\boldsymbol a}$ is a compact
connected Abelian group and has dimension 1. Moreover
$\mathbf{\Sigma}_{\boldsymbol a}^{*}\cong H_{\boldsymbol a}$,
where
$$H_{\boldsymbol a}=\{\frac{m}{a_0a_1\cdots a_n}:n=0,1,\ldots;m\in\mathbb{Z}\},$$
is subgroup of $\mathbb{Q}$. Denote by $\mathcal{P}$ the set of
prime numbers.

\section{Lemmas}

Let $X$ be a locally compact Abelian group. Put $Y=X^*$,
$\tilde{\alpha}_{ij}\in Aut(Y),i,j=1,2,\ldots,n$. Let $f_i(y)$ be
some functions on $Y$. Recall that the
Skitovich--Darmois equation is an equation of the form:
\begin{equation}\label{S-D obsh f}
\prod_{i=1}^{n}f_i\left(\sum_{j=1}^{n}\tilde{\alpha}_{ij}u_j
\right)=\prod_{i=1}^{n}\prod_{j=1}^{n}f_i(\tilde{\alpha}_{ij}u_j),\quad
u_j \in Y.
\end{equation}
The proof of the main theorem is reduced to the studying of the
solutions of this equation. In order to prove the main result we
need some lemmas.

\begin{lem}\label{Lemma fe}(\cite{Mazur}).
Let $X$ be a second countable locally compact Abelian group,
$\xi_i, i=1,2,\ldots,n,$ be independent random variables with
values in $X$, and distributions $\mu_i$. The linear forms
$L_j=\sum_{i=1}^{n}\alpha_{ij}\xi_i,j=1,2,\ldots,n,$ where
$\alpha_{ij}\in Aut(X)$, are independent if and only if the
characteristic functions $\hat{\mu}_i(y),i=1,2,\ldots,n,$ satisfy
equation (\ref{S-D obsh}), which takes the form
\begin{equation}\label{S-D obsh}
\prod_{i=1}^{n}\hat{\mu}_i\left(\sum_{j=1}^{n}\tilde{\alpha}_{ij}u_j\right)=
\prod_{i=1}^{n}\prod_{j=1}^{n}\hat{\mu}_i(\tilde{\alpha}_{ij}u_j),\quad
u_j \in Y.
\end{equation}

\end{lem}

\begin{lem}(\cite{Mazur}).\label{Finite 0}
Let $X$ be a direct product of groups $\mathbb{Z}(p^{k_p})$, where
$k_p\geq 0$, i.e. $X=\mathbf{P}_{p\in
\mathcal{P}}\mathbb{Z}(p^{k_p})$. Let $\xi_i, i=1,2,\ldots,n,$ be
independent random variables with values in $X$  and distributions
$\mu_i$. Then the independence of the linear forms
  $L_j=\sum_{i=1}^{n}\alpha_{ij}\xi_i,$
where  $\alpha_{ij}\in
Aut(X),\alpha_{1j}=\alpha_{i1}=I,i,j=1,2,\ldots,n$ implies that
$\mu_i=E_{x_i}\ast m_K$, where $K$ is a compact subgroup of $X$,
$x_i\in X$, $i=1,2,\ldots,n$.
\end{lem}

Taking into the account that $X=\mathbf{P}_{p\in
\mathcal{P}}\mathbb{Z}(p^{k_p})$ if and only if $Y$ is a weak
direct product of the groups $\mathbb{Z}(p^{k_p})$, where $k_p\geq
0$, i.e. $Y=\mathbf{P}^{*}_{p\in \mathcal{P}}\mathbb{Z}(p^{k_p})$,
by lemmas \ref{Lemma fe} and \ref{Finite 0} we obtain

\begin{cor}\label{Finite}
Let $Y$ be a discrete Abelian group of the form
$Y=\mathbf{P}^{*}_{p\in \mathcal{P}}\mathbb{Z}(p^{k_p})$, where
$k_p\geq 0$. Let $\hat{\mu}_i(y),i=1,2,\ldots,n$, $n\geq 2,$ be
the characteristic functions on $Y$, satisfying equation (\ref{S-D
obsh}), where $\tilde{\alpha}_{ij}\in Aut(Y),\tilde{\alpha}_{1j}=
\tilde{\alpha}_{i1}=I,i,j=1,2,\ldots,n$. Then
$\hat{\mu}_i(y)=(x_i,y)\hat{m}_K(y),y\in Y,$ where $K$ is a
compact subgroup of $X$, $x_i\in X$, $i=1,2,\ldots,n$.
\end{cor}

The following lemma states that an analogue of the
Skitovich--Darmois theorem for three linear forms of three
independent random variables  holds on the circle group if we
assume that the characteristic functions of the random variables
do not vanish.

\begin{lem}(\cite{MF})\label{MF}
Assume that $X=\mathbb{T},\alpha_{ij}\in Aut(X),i,j=1,2,3$. Let
$\xi_i, i=1,2,3,$ be independent random variables with values in
$X$ and  distributions $\mu_i$, such that their characteristic
functions do not vanish. Suppose that
$L_j=\sum_{i=1}^{3}\alpha_{ij}\xi_i,j=1,2,3,$ are independent.
Then $\mu_i=E_{x_i},$ $x_i\in X,$ $i=1,2,3$.
\end{lem}

By Lemmas \ref{Lemma fe} and \ref{MF} we obtain

\begin{cor}\label{3 formi}
Assume that $Y=\mathbb{Z}$. Let $\hat{\mu}_i(y),i=1,2,3$, $n\geq
2,$ be non-vanishing characteristic functions on $Y$ satisfying
the equation

$$\hat{\mu}_1(u_1+u_2+u_3)\hat{\mu}_2(u_1-u_2-u_3)\hat{\mu}_3(u_1+u_2-u_3)=$$
\begin{equation}\label{S-D T}
=\hat{\mu}_1(u_1)\hat{\mu}_1(u_2)\hat{\mu}_1(u_3)\hat{\mu}_2(u_1)\hat{\mu}_2(-u_2)\hat{\mu}_2(-u_3)
\hat{\mu}_3(u_1)\hat{\mu}_3(u_2)\hat{\mu}_3(-u_3), $$$$u_i \in
Y,i=1,2,3.
\end{equation}
Then $\hat{\mu}_i(y)=(x_i,y),x_i\in X,i=1,2,3,y\in Y$.
\end{cor}

\begin{lem}(\cite[Lemma $13.20$]{F book})\label{Lemma compGr}
Let $X$ be a second countable compact Abelian group. Suppose that there exists an
automorphism $\delta \in Aut(X)$ and an element $\tilde{y} \in Y$,
such that the following conditions are satisfied:

i)$Ker(I-\tilde{\delta})=\{0\}$;

ii)$(I-\tilde{\delta})Y\cap \{0;\pm \tilde{y},\pm 2\tilde{y}\}=\{0\}$;

iii) $\tilde{\delta}\tilde{y}\neq-\tilde{y}$.

Then for all $n\geq2$ there exist independent identically
distributed random variables $\xi_i,i=1,2,\ldots,n,$ with values
in $X$ and  distribution $\mu \not\in I(X)*\Gamma(X)$, such that
the linear forms
$L_j=\xi_1+\sum_{i=2}^{n}\delta_{ij}\xi_i,j=1,2,\ldots,n,$ where
$\delta_{ij}=I,i\neq j,\delta_{ii}=\delta,$ are independent.
\end{lem}

It is convenient for us to formulate as a lemma the following
simple statement.
\begin{lem}\label{PosDefFunc}
Let $Y$ be a second countable  discrete Abelian group, $H$ be a subgroup of $Y$,
$f(y)$ be a function on $Y$ of the form
\begin{equation}\label{pdfunc}
f(y)=
\begin{cases}
1,\quad y\in H; \\
c,\quad y \not\in H,
\end{cases}
\end{equation}
where $0<c<1$. Then  $f(y)$ is a positive defined  function.

\end{lem}
\textbf{Proof.} Consider the distribution $\mu=cE_{0}+(1-c)m_{G}$
on the group $X$, where $G=A(X,H)$. It is easy to see that
$f(y)=\hat{\mu}(y)$. Hence, $f(y)$ is a positive definite
function. $\blacksquare$.

The following lemma for $n=2$ was proved in \cite{FG 2000}.

\begin{lem}\label{Theorem CompConnectGroup}
Let $X$ be a second countable  compact connected Abelian group, such that $f_2\in
Aut(X)$. Then there exist independent random variables
$\xi_i,i=1,2,\ldots,n,$ with values in $X$ and distributions $\mu_i\not\in I(X)\ast
\Gamma(X)$, and automorphisms $\alpha_{ij}\in Aut(X)$, such that
the linear forms
$L_j=\sum_{i=1}^{n}\alpha_{ij}\xi_i,j=1,2,\ldots,n,$ are
independent.
\end{lem}
\textbf{Proof.} Two cases are possible: 1. $f_p \in Aut(X)$   for
all prime number $p$; 2. $f_p\not\in Aut(X)$ for a prime number
$p$.

\textbf{1.} Consider the first case. It is well-known that if $X$
is a compact Abelian group $X$ such that $f_p\in Aut(X)$ for all
prime $p$, then
\begin{equation}\label{CCG p}
X\cong (\mathbf{\Sigma}_{\boldsymbol a})^{\mathfrak{n}},
\end{equation}
where ${\boldsymbol a}=(2,3,4,\ldots)$,\quad (\cite[(25.8)]{H R}).
It is obvious that it suffices to prove the lemma for the group of
the form $X=\mathbf{\Sigma}_{\boldsymbol a},{\boldsymbol
a}=(2,3,4,\ldots)$. Then the group $Y$ is topologically isomorphic
to the group $\mathbb{Q}$. Let $p$ and $q$ be different prime
numbers. Let $H$ be a subgroup of $Y$ of the form
$H=\{\frac{m}{q^k}\}_{m,k\in\mathbb{Z}}$. Put
$G=H^*,K=A(G,H^{(p)})$. Since numbers $p$ and $q$ are relatively
prime, it follows that $H\neq H^{(p)}$. Consider on the  group $H$
the function
\begin{equation}
f(y)=
\begin{cases}
1,\quad y\in H^{(p)},\\
c,\quad y\not\in H^{(p)}, \label{Def lambda}
\end{cases}
\end{equation}
where $0<c<1$. By   Lemma \ref{PosDefFunc}  $f(y)$ is  a
 positive definite function.

Consider  on the group  $Y$ the function
\begin{equation}
g(y)=
\begin{cases}
f(y),\quad y\in H,\\
0,\quad\quad \ y\not\in H.
\end{cases}
\end{equation}
The function $g(y)$ is a positive definite function (\cite[Theorem
2.12]{F book}). By the Bohner theorem there exists a distribution
$\mu\in M^1(X)$ such that $\hat{\mu}(u)=g(y)$. It is obvious that
$\mu\not\in I(X)\ast \Gamma(X)$.

Let $\xi_i$ be independent random variables with values in $X$ and
distribution $\mu$. Put $s=p^2+q$. From the conditions of the
lemma it follows that $s\in Aut(X)$. Let us show that the linear
forms
$$L_1=\xi_1+p\xi_2+p\xi_3+\cdots+p\xi_n$$
$$L_2=p\xi_1+s\xi_2+p^2\xi_3+\cdots+p^2\xi_n$$
$$L_3=p\xi_1+p^2\xi_2+s\xi_3+\cdots+p^2\xi_n$$
$$\cdots\cdots\cdots\cdots\cdots\cdots\cdots\cdots\cdots\cdots\cdots $$
$$L_n=p\xi_1+p^2\xi_2+p^2\xi_3+\cdots+s\xi_n$$
are independent. By the Lemma \ref{Lemma fe} it suffices to show that the following equation holds:
$$\hat{\mu}(u_1+pu_2+pu_3+\cdots+pu_n)\hat{\mu}(pu_1+su_2+p^2u_3+\cdots+p^2u_n)\times\cdots$$
\begin{equation}
\times\hat{\mu}(pu_1+p^2u_2+\cdots+su_n)=\hat{\mu}(u_1)\hat{\mu}(pu_2)\hat{\mu}(pu_3)
\cdots\hat{\mu}(su_n). \label{S-D dlya f}
\end{equation}

From (\ref{Def lambda}) it follows that
\begin{equation}
g(y+pt)=g(y),\quad y,t\in H. \label{g(y+pt)}
\end{equation}
Using (\ref{g(y+pt)}) it is easy to show that if $u_i\in H$, then equation (\ref{S-D dlya f})
becomes a equality. So it suffices to consider the case when $u_i\not\in H$ for some $i$.
It is easy to see that in this case the right-hand side of equation (\ref{S-D dlya f}) vanishes.

Let us show that the left-hand side of equation (\ref{S-D dlya f}) vanishes too.
Assume the converse, i.e. that the left-hand side of equation (\ref{S-D dlya f}) does not vanish.
Then the following system of equations holds:
\begin{equation}
\begin{cases}
u_1+pu_2+pu_3+\ldots+pu_n=h_1,\\
pu_1+su_2+p^2u_3+\ldots+p^2u_n=h_2,\\
\ldots\ldots\ldots\ldots\ldots\ldots\ldots\ldots\ldots\ldots\ldots \\
pu_1+p^2u_2+p^2u_3+\ldots+su_n=h_n, \label{sys n}
\end{cases}
\end{equation}
where $h_i\in H$.

Add the first equation of the system (\ref{sys n}) multiplied by $(-p)$
to the each equation of the system (\ref{sys n}) starting from the second.
We obtain that $qu_i=h_i-ph_1,i=2,3,\ldots,n$. Thus $u_i\in H,i=2,3,\ldots,n$.
From this and from the first equation of the system (\ref{sys n}) it follows that that $u_1\in H$.
Finally we obtain that $u_i\in H,i=1,2,\ldots,n.$ This contradicts the assumption.

\textbf{2.} Assume that
\begin{equation}
f_p\not\in Aut(X) \label{fp not autt}
\end{equation}
for some prime number $p$. Suppose that $p$ is the smallest one
satisfying condition (\ref{fp not autt}). Since $X$ is a connected
group, we have $X^{(n)}=X$ for all natural $n$. Hence if
$f_p\not\in Aut(X)$, then $Kerf_p\neq \{0\}$.

From the condition of the lemma it follows that $p\geq 3$. Put
$a=1-p$. Since $p$ is a smallest natural number satisfying
condition (\ref{fp not autt}), we obtain $f_{-a}\in Aut(X)$. Hence
$f_{a}\in Aut(X)$. Note that $Ker f_p=A(X,Y^{(p)})$. It implies
that $Y^{(p)}\neq Y$. Let $\tilde{y}\in Y^{(p)}$ and verify that
the automorphism $\delta=f_a$ and the element $\tilde{y}$ satisfy
to conditions of Lemma \ref{Lemma compGr}. We have
$\tilde{f}_a=f_a$ and $I-\tilde{f}_a=\tilde{f}_p$. Since $Y$ is
torsion-free group, it follows that $Ker(I-\tilde{f}_{a})=\{0\}$,
i.e. condition (i) holds. Thus $(I-\tilde{f}_a)Y=Y^{(p)}$. From
$p\geq 3$ it follows that numbers $2$ and $p$ are relatively
prime. Hence there are integers $m$ and $n$ such that $2m+pn=1$.
Thus $y=2my+pny$. So if $\tilde{y}\not\in Y^{(p)}$, then
$2\tilde{y}\not\in Y^{(p)}$ too. It implies that condition (ii)
holds. Taking in the account that $Y$ is torsion-free group, it is
obvious that condition (iii) holds. We use Lemma \ref{Lemma
compGr} and obtain the assertion of the lemma.$\blacksquare$

\section{Main theorem}

\begin{thm}\label{Theorem Sol p}
Let $X=\mathbf{\Sigma}_{\boldsymbol a}$ be an ${\boldsymbol
a}$-adic solenoid.

$\mathbf{1}.$ Assume that $f_p\not\in Aut(X)$ for all prime numbers $p$.
Let $\xi_i, i=1,2,3,$ be independent random variables with values
in $X$ and distributions $\mu_i$. Then the independence of the
linear forms $L_j=\sum_{i=1}^{3}\alpha_{ij}\xi_i,$ where
$\alpha_{ij}\in Aut(X),i,j=1,2,3$, implies that at least one
distribution $\mu_{i}\in I(X)$.

$\mathbf{2}.$ Assume that $f_p\in Aut(X)$ for a prime number $p$. Then
there are independent random variables $\xi_i,i=1,2,3,$ with
values in $X$ and distributions $\mu_i\not\in \Gamma(X)\ast I(X)$,
and automorphisms $\alpha_{ij}\in Aut(X)$, such that the linear
forms $L_j=\sum_{i=1}^{3}\alpha_{ij}\xi_i,j=1,2,3,$ are
independent.
\end{thm}

It should be noted that an example of a group such that
$f_p\not\in Aut(X)$ for all prime number $p$ is the group
$\mathbf{\Sigma}_{\boldsymbol a},{\boldsymbol
a}=(2,3,5,7,\ldots)$. Its character group
$\mathbf{\Sigma}_{\boldsymbol a}^{*}\cong \{\frac{m}{p_1p_2\cdots
p_k} \ :m\in\mathbb{Z}, \ p_j \ \mbox{are different prime
numbers}\}$.

An example of a group such that $f_p\in Aut(X)$ for a prime number
$p$ is the group $\mathbf{\Sigma}_{\boldsymbol a},{\boldsymbol
a}=(2,2,2\ldots)$. Its character group
$\mathbf{\Sigma}_{\boldsymbol a}^{*}\cong \{\frac{m}{2^k} :m,
k\in\mathbb{Z}\}$.

The proof of Theorem 4.1 is divided into two parts. In the first
part we use Corollaries  \ref{Finite} and \ref{3 formi}. In the
second part we use Lemma \ref{Theorem CompConnectGroup}.

\textbf{Proof.} \textbf{1}.  Suppose that $f_p\not\in Aut(X)$ for all prime
numbers $p$. This implies that $Aut(X)=\{I,-I\}$. It is easy to
show that the case of arbitrary  linear forms $L_j$ is reduced to
the case when $L_j$ are of the form
$$L_1=\xi_1+\xi_2+\xi_3,$$
\begin{equation}
L_2=\xi_1-\xi_2+\xi_3, \label{formiT}
\end{equation}
$$L_3=\xi_1-\xi_2-\xi_3.$$

Note that $Y$ is topologically isomorphic to a subgroup of
$\mathbb{Q}$. To avoid introducing new notation  we will suppose
that $Y$ is a subgroup of $\mathbb{Q}$. By Lemma \ref{Lemma fe}
the independence of the linear forms (\ref{formiT}) implies that
equation (\ref{S-D T}), where $Y$ is a subgroup of $\mathbb{Q}$,
holds.

Note that, since $f_2\not\in Aut(X)$, we have that the partition
of $Y$ into the cosets of $Y^{(2)}$ consists of two cosets:
$Y^{(2)}$ and $\tilde{y}+Y^{(2)}$, where $\tilde{y} \not\in
Y^{(2)}$.

Put $N_{i}=\{y\in Y: \hat{\mu}_i(y)\neq0\}, N=\cap_{i=1}^3 N_i$.
We infer from (\ref{S-D T}) that $N$ is a subgroup in $Y$.
Moreover, it is easy to see from (\ref{S-D T}), that $N$ has a
property:
\begin{equation}\label{proph}
\mbox{if } 2y\in N, \mbox{ then } y\in N.
\end{equation}
There are two cases: $N\neq\{0\}$ and $N=\{0\}$.

\textbf{A.} Assume that $N\neq\{0\}$. Suppose that $t_1$ and $t_2$ belong to the same coset of
$Y^{(2)}$ in $Y$. Then there exist $\hat{u}_1$ and $\hat{u}_2$,
such that $\hat{u}_1+\hat{u}_2=t_1,\hat{u}_1-\hat{u}_2=t_2$. Putting first $u_1=\hat{u}_1,u_2=\hat{u}_2,u_3=0$
in (\ref{S-D T}), then $u_1=\hat{u}_1,u_2=-\hat{u}_2,u_3=0$ in (\ref{S-D T}), and equating the right-hand sides of obtained equations, we infer:
$$|\hat{\mu}_1(t_1)||\hat{\mu}_2(t_2)||\hat{\mu}_3(t_1)|=|\hat{\mu}_1(t_2)||\hat{\mu}_2(t_1)||\hat{\mu}_3(t_2)|.$$
Reasoning the same way, it is easy to see that if $t_1$ and $t_2$ belong to the same coset of
$Y^{(2)}$ in $Y$, then the following equation holds:
\begin{equation}
|\hat{\mu}_{i_1}(t_1)||\hat{\mu}_{i_2}(t_2)||\hat{\mu}_{i_3}(t_2)|=|\hat{\mu}_{i_1}(t_2)|
|\hat{\mu}_{i_2}(t_1)||\hat{\mu}_{i_3}(t_1)|, \label{perestanovka
argumentov}
\end{equation}
where all $i_j$ are pairwise different.

Put $\nu_i=\mu_i\ast\bar{\mu}_i,i=1,2,\ldots,n$. Then
$\hat{\nu_i}(y)=|\hat{\mu}_i(y)|^2, y\in Y$. Functions
$\hat{\nu_i}(y)$ are nonnegative and also satisfy equation
(\ref{S-D T}). It suffices to show that $\hat{\nu}_i(y)$ are
characteristic functions of the idempotent distributions. This
implies that $\hat{\mu}_i(y)$   are also characteristic functions
of the idempotent distributions.

Now we will show that $N_i=N,i=1,2,3$. Assume the converse. Then
there exists $y_1\in N_{i_1}$ such that either $y_1\not\in
N_{i_2}$ or $y_1\not\in N_{i_3}$, where all $i_j$ are pairwise
different. Put $t_1=y_1,t_2=y_2$, where $y_2\in N$ and $y_1$,$y_2$
belong to the same coset of the $Y^{(2)}$ in $Y$, in
(\ref{perestanovka argumentov}). We can make such choice. Indeed,
on the one hand $N\cap Y^{(2)}\neq\{0\}$ because $N$ is a subgroup
and $N\neq\{0\}$ by the assumption. On the other hand there exists
$y\neq 0$ such that $y\in N\cap (\tilde{y}+ Y^{(2)})$. Indeed, if
$N\subset Y^{(2)}$, then taking in the account (\ref{proph}) we
infer that there exists $y^{\prime}\in Y$ such that
$y=2^{k}y^{\prime}$. This contradicts to the fact that there are
no $y\in Y$ such that $y$ is infinitely divisible by $2$. We infer
that the left-hand side of equation (\ref{perestanovka
argumentov}) is equal to a positive number, and the right-hand side
of equation (\ref{perestanovka argumentov}) is equal to zero. This
is a contradiction. So we have that $N_i=N,i=1,2,3$.

Note that if $y\in N$, then $\hat{\nu}_i(y)=1,$ $i=1,2,3.$
Indeed, let $y_0\in N$. Consider the subgroup $H$ of $Y$ generated
by $y_0$. Note that $H\cong \mathbb{Z}$. Consider the restriction
of equation (\ref{S-D T}) to the subgroup $H$. Using Corollary
\ref{3 formi} we obtain that $\hat{\nu}_i(y)=1,i=1,2,3,y\in H$.

Taking in the account that the characteristic functions
$\hat{\nu}_{i}(y)$ are $N$-invariant, consider the equation
induced by  equation (\ref{S-D T}) on the factor-group $Y/N$. Put
$f_i([y])=\hat{\nu}_{i}([y])$. Note that if $H$ is an arbitrary nontrivial subgroup of $Y$, then $Y/H$ is topologically isomorphic to a group of the form $\mathbf{P}^{*}_{p\in
\mathcal{P}}\mathbb{Z}(p^{k_p})$, where $k_p\geq 0$. In particular, this holds for the factor-group $Y/N$. Hence, by
Corollary \ref{Finite} we infer that $f_i([y])$ are characteristic
functions of some idempotent distributions. This implies that all
distributions $\mu_i$ are idempotent.

\textbf{B.} Consider the case $N=\{0\}$.

Put first $u_2=0,u_3=u_1=y$, after $u_3=0,u_1=u_2=y$, and finally
$u_1=0,u_2=u_3=y$ in (\ref{S-D T}), we infer respectively:
\begin{equation}
\hat{\mu}_1(2y)=\hat{\mu}_1^2(y)|\hat{\mu}_2(y)|^2|\hat{\mu}_3(y)|^2,\quad y\in Y.\label{2y=yyy 1}
\end{equation}

\begin{equation}
\hat{\mu}_2(2y)=|\hat{\mu}_1(y)|^2\hat{\mu}_2^2(y)|\hat{\mu}_3(y)|^2,\quad y\in Y.\label{2y=yyy 2}
\end{equation}
\begin{equation}
\hat{\mu}_3(2y)=|\hat{\mu}_1(y)|^2|\hat{\mu}_2(y)|^2\hat{\mu}_3^2(y),\quad y\in Y.\label{2y=yyy 3}
\end{equation}

Note that
\begin{equation}
\hat{\mu}_i(2y)=0,y\in Y,y\neq0,i=1,2,3. \label{mu(2y)=0}
\end{equation}
Indeed, if $\hat{\mu}_{i_0}(2y_0)\neq0$ for some $y_0\in Y,y_0\neq
0,$ and $i_0$, then from equalities (\ref{2y=yyy 1})-(\ref{2y=yyy
3}) it follows that $\hat{\mu}_i(y_0)\neq0, i=1,2,3$. This
contradicts to $N\neq \{0\}$.

Show that at least one distribution $\mu_i=m_X$. Assume the converse.
Then there exist $t_1\neq 0,t_2\neq 0,t_3\neq 0$, such that
\begin{equation}\label{pst0}
\hat{\mu}_1(\pm t_1)\hat{\mu}_2(\pm t_2)\hat{\mu}_3(\pm t_3)\neq 0.
\end{equation}
From equality (\ref{mu(2y)=0}) it follows that $t_i\in
\tilde{y}+Y^{(2)}$. From $N=\{0\}$ it follows that $\pm
t_i,i=1,2,3,$ do not coincide. Without lost of generality assume
that $t_1\neq\pm t_2$. Note that for all elements
$y^{\prime},y^{\prime\prime} \in \tilde{y}+Y^{(2)}$ we have
$y^{\prime}+y^{\prime\prime} \in Y^{(2)}$. Moreover, for any two
elements $y^{\prime},y^{\prime\prime}\in \tilde{y}+Y^{(2)}$ there
are two possibilities: either $y^{\prime}+y^{\prime\prime}\in
Y^{(4)}$ or $y^{\prime}-y^{\prime\prime}\in Y^{(4)}$.

Put $y_i=t_i,i=1,2,3,$ if $t_1+t_2\in Y^{(4)},$ and put $y_1=t_1,y_2=-t_2,y_3=t_3,$ if $t_1-t_2\in Y^{(4)}$.
Note that $y_1+y_2\in Y^{(4)},y_1+y_2\neq0$.
For an element $y_0\in Y^{(2)}$ denote by $\frac{y_0}{2}$ such element of $Y$,
that $2\frac{y_0}{2}=y_0$. Thus we infer that $\frac{y_1+y_2}{2}\in Y^{(2)},\frac{y_1+y_2}{2}\neq 0$.

Consider the system of equations
\begin{equation}\label{pst1}
\begin{cases}
u_1+u_2+u_3=y_1,\\u_1-u_2-u_3=y_2,\\u_1+u_2-u_3=y_3.
\end{cases}
\end{equation}
Taking in the account that $y_i\in \tilde{y}+Y^{(2)},i=1,2,3,$ it is easy to see that the system of equations (\ref{pst1}) has the following solutions
\begin{equation}\label{pst2}
\begin{cases}
u_1=\frac{y_1+y_2}{2},\\u_2=\frac{y_3-y_2}{2},\\u_3=\frac{y_1-y_3}{2}
\end{cases}
\end{equation}
Put the solutions of (\ref{pst2}) in equation (\ref{S-D T}).
Taking into the account (\ref{pst0}) we infer that the right-hand
side of (\ref{S-D T}) is not equal to 0. This implies that
\begin{equation}\label{pst3}
\hat{\mu}_1 \left(\frac{y_1+y_2}{2}\right)\hat{\mu}_2\left(\frac{y_3-y_2}{2}\right)
\hat{\mu}_3\left(\frac{y_1-y_3}{2}\right)\neq 0.
\end{equation}

It follows from inequality (\ref{pst3})   that
$\mu_{1}(\frac{y_1+y_2}{2})\neq 0$. However, we have $\frac{y_1+y_2}{2}\in
Y^{(2)}$. This contradicts to (\ref{mu(2y)=0}).

Note that we proved more: if $N\neq \{0\}$ then all distributions
$\mu_i$ are idempotent, and if $N= \{0\}$ at least one
distribution $\mu_i$ is the Haar distribution on $X$.

\textbf{2}. Now consider the case $f_p\in Aut(X)$ for some prime $p$.
If $f_2\in Aut(X)$, then the statement follows from the Lemma \ref{Theorem CompConnectGroup}.
Assume that $f_2\not\in Aut(X)$.
Then two cases are possible: $p-1=4k$ and $p+1=4k$. Let us study the first case.

Consider the function $\rho(x)$ on $X$ defined by the equation
$$\rho(x)=1+Re(x,y_0),$$
where $y_0\in Y, y_0\neq 0$. It is obvious that $\rho(x)\geq
0,x\in X,$ and $\int_X \rho(x)dm_X(x)=1$. Let $\mu$ be a
distribution on $X$ with the density $\rho(x)$ with respect to
$m_X$. It is obvious that $\mu\not\in \Gamma(X)\ast I(X)$. The
characteristic function of the distribution $\mu$ is of the form:
\begin{equation}\label{def mu}
\hat{\mu}(y)=
\begin{cases}
1,& y=0,\\
\frac{1}{2},& y=\pm y_0,\\
0,& y\not\in \{0,y_0,-y_0\}.
\end{cases}
\end{equation}

Let $\xi_i,i=1,2,3,$ be independent identically distributed random variables with values in $X$ and distribution $\mu$.
Let us verify that the linear forms $L_1=\xi_1+\xi_2+\xi_3,L_2=\xi_1+p\xi_2+\xi_3,L_3=\xi_1+\xi_2+p\xi_3$
 are independent.
By Lemma \ref{Lemma fe} it suffices to prove that $\hat{\mu}(y)$ satisfies equation (\ref{S-D obsh}), which takes the form
\begin{equation}
\hat{\mu}(u+v+t)\hat{\mu}(u+pv+t)\hat{\mu}(u+v+pt)=\hat{\mu}^3(u)\hat{\mu}^2(v)\hat{\mu}^2(t)\hat{\mu}(pv)\hat{\mu}(pt),\label{S-D mu}
\end{equation}
where $u,v,t \in Y$. We will show that equation (\ref{S-D mu})
holds. It is obvious, that it suffices to consider the case, when at
least two of three elements $u,v,t$ are not equal to $0$. It is
easy to see that in this case the right-hand side of equation
(\ref{S-D mu}) is equal to 0. Let us show that the left-hand side
of equation (\ref{S-D mu}) vanishes too.

Suppose that there are some elements $u,v,t$ such that the left-hand side of equation (\ref{S-D mu}) does not vanish. Then there exist some $h_i\in \{0,y_0,-y_0\},i=1,2,3,$ such that $u,v,t$ satisfy the system of equations
\begin{equation}\label{syst 1,p}
\begin{cases}
u+v+t=h_1,\\
u+pv+t=h_2,\\
u+v+pt=h_3.
\end{cases}
\end{equation}

It is easy to obtain from (\ref{syst 1,p}) that
\begin{equation}
(p-1)v,(p-1)t\in \{0,\pm y_0,\pm 2y_0\}. \label{y_0 in}
\end{equation}

Relationship (\ref{y_0 in}) fails because of $(p-1)=4k$, but $y_0 \not\in Y^{(2)}$. From this it follows that the left-hand side of equation (\ref{S-D mu}) is equal to $0$.

The second case can be studied similarly. In this case we have to consider the linear forms $L_1=\xi_1+\xi_2+\xi_3,L_2=\xi_1-p\xi_2+\xi_3,L_3=\xi_1+\xi_2-p\xi_3$.

The theorem is completely proved.

$\blacksquare$

This research was conducted in the frame of the project "Ukrainian
branch of the French-Russian Poncelet laboratory" - "Probability
problems on groups and spectral theory" and was supported in part by the grant ANR-09-BLAN-0084-01 and by the scholarship of the French Embassy in Kiev.


\begin{thebibliography}{99}




\bibitem{F 1992} G. M. Feldman, On the Skitovich-Darmois theorem for finite abelian groups.// \textit{Theory Probabl. Appl.} \textbf{37} (1992), 621-631.


\bibitem{FG 2000} G. M. Feldman and P. Graczyk, On the Skitovich-Darmois theorem on compact Abelian groups.// \textit{J. of Theoretical Probability,} \textbf{13} (2000), 859-869.

\bibitem{F 2005} On a characterization theorem for locally compact abelian
groups.// \emph{Probab. Theory Relat. Fields,} Vol. 133, (2005),
345–357.

\bibitem{FG 2005} G. M. Feldman and P. Graczyk, On the Skitovich-Darmois theorem for discrete Abelian groups. //\textit{Theory Probab. Appl.} \textbf{49} (2005), 527-531.

\bibitem{FM 2008} G. Feldman, M. Myronyuk , The Skitovich-Darmous equation in the classes of continuous and measurable functions.// Aequationes Mathematicae 75 (2008), 75-92.




\bibitem{F book} Feldman G., \textit{Functional equations and Characterizations problems on locally compact Abelian Groups} // EMS Tracts in Mathematics Vol. 5. Zurich. European Mathematical Society (EMS) 2008. -- 268 pp.

\bibitem{MF} G.M.Feldman,  M.V. Myronyuk., Independent linear statistics on the
cylinders. // arXiv:1212.2772v1, (2012).






\bibitem{KLR} A. M. Kagan, Yu. V. Linnik, C. Radhakrishna Rao,\emph{ Characterization Problems in Mathematical Statistics.}  //
1973. --New York : Wiley, - Wiley series in probability and mathematical statistics.


\bibitem{Mazur Fin} I.P. Mazur, The Skitovich-Darmois theorem for finite Abelian groups (the case of n linear forms of n random variables)// \emph{U.M.J.}  2011.-- \textbf{4}. -- 1512-1523.

\bibitem{Mazur} I.P. Mazur, The Skitovich-Darmois theorem for discrete and compact totally disconnected Abelian
groups.// arXiv:1112.1488v1, (2011).





\bibitem{H R} Hewitt E. and A.Ross, \emph{Abstract harmonic analysis.} Vol.1 (Springer-Verlag, Berlin, Gottingen, Heildelberg, 1963).




\bibitem{Parth} K. R. Parthasarathy, \emph{Probability Measures on Metric Spaces.} Chelsea Publishing (AMS), 1967.
\end{thebibliography}
\end{document}